\documentclass[12pt]{amsart}

\usepackage{amsfonts}
\usepackage{amssymb}
\usepackage{amsmath,amsthm}
\usepackage[mathscr]{eucal}
\usepackage[latin1]{inputenc}
\usepackage[notcite,notref]{%showkeys
}

\numberwithin{equation}{section}

\theoremstyle{plain}
\newtheorem{teo}{Theorem}
\newtheorem{lem}{Lemma}

\theoremstyle{remark}
\newtheorem{rem}{Remark}

\theoremstyle{definition}
\newtheorem{defi}{Definition}

\title{Embeddings for anisotropic Besov spaces}

\author{F.J. P\'erez L\'azaro}\thanks{Thanks to V.I. Kolyada for his useful ideas and suggestions.}
\thanks{Research supported in part by grant MTM2006-13000-C03-03 of the DGI, Spain}
\address{Departamento de Matem\'aticas y Computaci\'on,
Universidad de La Rioja, Edificio J.~L.~Vives, Calle Luis de Ulloa
s/n, 26004 Logro\~no, Spain}
\email{francisco-javier.per@dmc.unirioja.es}

\date{}

\begin{document}
%%%%%%%%%%%%%%%%%%%%%%%%%%%%%%%%%%%%%%%%%%%%%%%%%%

\maketitle

%%%%%%%%%%%%%%%%%%%%%%%%%%%%%%%%%%%%%%%%%%%%%%%%%%

\begin{abstract}
  We prove embedding theorems for fully anisotropic Besov spaces.
  More concrete, inequalities between modulus of continuity in different metrics  and
  of Sobolev type are obtained. Our goal is to get sharp estimates for some anisotropic
  cases previously unconsidered.
\end{abstract}

\section{Introduction}
This paper places in the theory of embeddings of spaces of
differentiable functions in several variables. Our objective is to
prove embeddings for anisotropic Besov spaces. The main result is a
sharp embedding of different metrics (which is a generalization of
the classical embedding between Nikol'ski\u{\i} classes
\cite{Nik51}) for Besov type spaces with all the parameters that can
be different for each variable. Embeddings for anisotropic Besov
spaces in Lorentz spaces are proved also. This work continues
\cite{KoPe}, where embeddings for anisotropic Sobolev spaces were
found.

In order to specify better the results, let us recall  an historical
preview.

 The study of spaces of
differentiable functions in several variables with fractional index
of smoothness was begun by Nikol'ski\u{\i}, connected with problems
in Approximation Theory. He obtained some embeddings for the classes
$H_p^{r_1,\ldots,r_n}$, characterized for H\"older conditions in
$L^p$ for the differences of the derivatives of various orders. In
particular, this analogous of the theorem of Hardy-Littlewood was
proved in \cite{Nik51}: Let $r_j>0$ ($j=1,\ldots,n$),
$r:=n(\sum_{j=1}^nr_j^{-1})^{-1}$, $1\le p <q <\infty$. If
$\varkappa:=1-\frac{n}{r}(\frac{1}{p}-\frac{1}{q})>0$ and
$\alpha_j:=\varkappa r_j$ ($j=1,\ldots,n$), then
\begin{equation}\label{metricash}
  H_p^{r_1,\ldots,r_n}(\mathbb{R}^n)\hookrightarrow
  H_q^{\alpha_1,\ldots,\alpha_n}(\mathbb{R}^n).
\end{equation}
Later, a theory of similar spaces was built by Besov: the scale of
the so called $B$-spaces, introduced by him. In \cite[vol.2,
pg.62]{BeIlNi}, using the previous notation and $1\le \theta\le
\infty$, it is obtained that
\begin{equation}\label{metricasb}
  B_{p,\theta}^{r_1,\ldots,r_n}(\mathbb{R}^n)\hookrightarrow
  B_{q,\theta}^{\alpha_1,\ldots,\alpha_n}(\mathbb{R}^n).
\end{equation}
 Observe that in \cite{Ko98},
Kolyada showed that (\ref{metricash}) can be proved using estimates
of non increasing rearrangements. Let us note that (\ref{metricash})
is the particular case $\theta=\infty$ in (\ref{metricasb}). In
fact, in \cite[vol.2, pg.62]{BeIlNi}, the embedding
(\ref{metricasb}) was proved in a more general form. It is
considered the case when, for each different variable $x_j$, is
taken a different metric $L^{p_j}$.

It is a logical continuation to consider, for each different
variable $x_j$, not only values $p_j$, but also values $\theta_j$
that can be different. However, the case when in the definition of
Besov space are included different values of the parameters
$\theta_j$ for each variable, had not been treated. Our objective is
to find the sharp parameters for the embedding in this case. The
main result in this paper is (see Theorem \ref{maint} below) an
inequality that implies the embedding
\begin{equation}\label{bminusabminus}
  b_{p_1,\ldots,p_n;\theta_1,\ldots,\theta_n}^{r_1,\ldots,r_n}(\mathbb{R}^n)\hookrightarrow
  b_{q_1,\ldots,q_n;\theta'_1,\ldots,\theta'_n}^{\alpha_1,\ldots,\alpha_n}(\mathbb{R}^n).
\end{equation}
We shall not specify here the conditions on the parameters. Let us
emphasize that the new in this theorem is to obtain the optimal
values $\theta'_j$, for the case when the parameters $\theta_j$ are
different.

As we shall see in Remark \ref{mejorlorentz}, the estimate obtained
in Theorem \ref{maint} is stronger than the embedding
(\ref{bminusabminus}).

In this paper we prove also Sovolev type inequalities. More
specifically, we obtain embeddings of Besov spaces into Lorentz
spaces. This kind of inequalities were proved and extended in
\cite{Ul,Pee,Go,He,BeIlNi,Ko98} and others. The results showed in
this paper (Theorem \ref{limit} for the embedding with limit
exponent and Theorem \ref{nolimit} for the embedding without limit
exponent) extend the previous results to the case of anisotropic
Besov spaces where the parameters $\theta'_j$ can be different.

Our methods are based on estimates of non increasing rearrangements.
The first works using this approach on the theory of embedding of
function classes are due to Ul'yanov at the end of sixties. Later,
these methods were mainly developed by Kolyada (see, for instance,
\cite{Ko93,Ko98,Ko01}). Here we use estimates of the non increasing
rearrangement of a function in terms of the modulus of continuity
obtained in \cite{Ko98}. In order to allow this approach work for
anisotropic Besov spaces we need to find a kind of sharp equilibrium
between the estimates. For this we apply the methods developed in
\cite{KoPe} and used also in \cite{Pe}.

This paper is organized as follows. Section \ref{secdos} contains
the basic definitions and notations. Section \ref{seclem}, the
lemmas used in the proofs. In section \ref{secteo}, the results of
this paper are presented, as well as their proofs and some remarks
about them.

\section{Definitions}\label{secdos}

Set $\mathbb{R}_+\equiv (0,+\infty)$. For $1\le p<\infty$ we denote
$\mathcal{L}^p\equiv L^p(\mathbb{R}_+,du/u)$; say also
$\mathcal{L}^\infty\equiv L^\infty(\mathbb{R}_+)$ (see \cite{He}).
Note that $L^\infty(\mathbb{R}_+)=L^\infty(\mathbb{R}_+,du/u)$.

Let $S_0(\mathbb{R}^n)$ be the class of all measurable and almost
everywhere finite functions $f$ on $\mathbb{R}^n$ such that for each
$y>0$,
\begin{equation*} \lambda_f (y) \equiv | \{x \in \mathbb{R}^n : |f(x)|>y \}| < \infty.  \end{equation*}

A non-increasing rearrangement of a function $f \in
S_0(\mathbb{R}^n)$ is a non-increasing function $f^*$ on
$\mathbb{R}_+$ that is equimeasurable with $|f|$. The rearrangement
$f^*$ can be defined by the equality
\begin{equation*} f^*(t) = \sup_{|E|=t} \inf_{x \in E} |f(x)| \,\quad,\quad 0<t<\infty \,. \end{equation*}

Assume that $0<q,p<\infty$ . A function $f \in S_0(\mathbb{R}^n)$
belongs to the Lorentz space $L^{q,p}(\mathbb{R}^n)$ if
\begin{equation*} \|f\|_{q,p} \equiv \left( \int_0^\infty \left( t^{1/q} f^*(t) \right)^p \frac{dt}{t} \right)^{1/p} < \infty \,. \end{equation*}
We have the inequality \cite[p.217]{BeSh}
\begin{equation*} \|f\|_{q,s} \le c \|f\|_{q,p} \quad (0<p<s<\infty), \end{equation*}
so that $L^{q,p} \subset L^{q,s}$ for $p<s$. In particular, for
$0<p\le q$
\begin{equation*} L^{q,p} \subset L^{q,q} \equiv L^q \,. \end{equation*}

The differences of degree $k$ in the direction of the variable $x_j$
are defined as
\begin{equation*}
  \Delta^k_j(h)f(x)\equiv \sum_{i=0}^k
  (-1)^{k-i}\binom{k}{i}f(x+ihe_j) \quad (h\in \mathbb{R}),
\end{equation*}
where $e_k$ is the unit coordinate vector. The modulus of continuity
in the metric $L^{q,p}$:
\begin{equation*}
  \omega_j^k(f;\delta)_{q,p}=\sup_{0<h\le\delta}\|\Delta^k_j(h)f\|_{q,p}.
\end{equation*}

As in \cite[pg.152,161]{Ni}, we define the Besov space in the
direction of the coordinate axe $x_j$.
\begin{defi}
  Let $r>0$, $1\le p<\infty$, $1\le \theta\le \infty$ and $1\le j\le
  n$. We define the space $B_{p,\theta;j}^r(\mathbb{R}^n)$ as the
  class of functions $f\in L^p(\mathbb{R}^n)$ such that
  \begin{equation*}
    \|f\|_{B_{p,\theta;j}^r}\equiv \|f\|_p+\|f\|_{b_{p,\theta;j}^r},
  \end{equation*}
  where the seminorm
  \begin{equation*}
    \|f\|_{b_{p,\theta;j}^r}\equiv \left(\int_0^\infty
    [h^{-r}\|\Delta_j^k(h)f\|_p]^\theta
    \frac{dh}{h}\right)^{1/\theta}=
    \left\|h^{-r}\|\Delta_j^k(h)f\|_p\right\|_{\mathcal{L}^\theta}.
  \end{equation*}
\end{defi}
It is well known that each election of an integer $k>r$ gives
equivalent seminorms. Furthermore, if we change the expression
$\|\Delta_j^k(h)f\|_p$ by the modulus of continuity
$\omega_j^k(f;h)_p$ we obtain equivalent seminorms also
(\cite[Chapter 4]{BeIlNi} and \cite[Chapter 4]{Ni}).

Moreover, the following inquality between seminorms holds (see, for
instance
  \cite[vol.2, pg. 64]{BeIlNi}). If $1\le \theta_1<\theta_2\le
  \infty$ then $\|.\|_{b_{p,\theta_2;j}^r}\le
  c\|.\|_{b_{p,\theta_1;j}^r}$. So, the bigger is $\theta$, the
  bigger is the corresponding Besov space.

Now we are going to present anisotropic Besov spaces.
\begin{defi}
  Let $n\in \mathbb{N}$, $r_j>0$, $1\le p_j<\infty$, $1\le
  \theta_j\le \infty$ $(j=1,\ldots,n)$. We say that $f\in
  b_{p_1,\ldots,p_n;\theta_1,\ldots,\theta_n}^{r_1,\ldots,r_n}(\mathbb{R}^n)$
  if $f\in S_0(\mathbb{R}^n)$ and the following seminorm is finite
  \begin{equation*}
    \sum_{j=1}^n\|f\|_{b_{p_j,\theta_j;j}^{r_j}}.
  \end{equation*}
\end{defi}
In the case $p_1=\cdots=p_n=p$ and $\theta_1=\cdots=\theta_n=\theta$
we use the notation $ b_{p;\theta}^{r_1,\ldots,r_n}(\mathbb{R}^n):=
b_{p,\ldots,p;\theta,\ldots,\theta}^{r_1,\ldots,r_n}(\mathbb{R}^n)$
and
$B_{p;\theta}^{r_1,\ldots,r_n}(\mathbb{R}^n):=L^p(\mathbb{R}^n)\cap
b_{p;\theta}^{r_1,\ldots,r_n}(\mathbb{R}^n)$. As usual
$H_{p}^{r_1,\ldots,r_n}(\mathbb{R}^n):=
B_{p;\infty}^{r_1,\ldots,r_n}(\mathbb{R}^n)$.
\section{Lemmas}\label{seclem}

The following lemma was proved in \cite[pg.167]{Ko98}. It presents
estimates of rearrangements in terms of the modulus of continuity.

%%%%%%%%%%%%%%%%%%%%%%%%%%%%%%%%%%%%%%%%%%%%%%%%%
% estimaciones del east journal
%%%%%%%%%%%%%%%%%%%%%%%%%%%%%%%%%%%%%%%%%%%%%%%%%
\begin{lem}\label{eastj}
Let $f$ be a locally integrable function in $S_0(\mathbb{R}^n)$,
$k_i \in \mathbb{N}$ and $p_i \in [1,\infty)$ $(i=1,\ldots,n)$.
Suppose that $\delta_i(t)$ $(1 \le i \le n)$ are positive functions
in $\mathbb{R}_+$ such that
\begin{equation}\label{l1}
  \prod_{i=1}^n \delta_i(t) =t \quad (t>0).
\end{equation}
Then, for all $0<t<s< \infty$,
\begin{equation}\label{l2}
  f^*(t) \le (2^k -1 ) f^*(s) + c \max_{i \in \{1,\ldots,n\}} t^{-1/p_i} \left( \frac{s}{t}\right)^{k_i} w_i^{k_i}(f; \delta_i(t))_{p_i},
\end{equation}
where $k=\max k_i$ and $c$ is a constant which only depends on $p_i$ and $k_i$.
\end{lem}

The previous lemma is formulated in \cite{Ko98} for functions in
$L^p(\mathbb{R}^n)$, and $p=p_1=\cdots=p_n$. If the reader checks it
carefully will realize that the proof is still valid without these
assumptions.

The aim of the next lemma is that given a function with some
properties of monotony and integrability we can majorized it for
another one with equivalent integrability properties and for which
its increasing and decreasing is controlled.

\begin{lem}\label{crec}
Let $\alpha >0$, $\theta \ge 1$. Let $\psi(t)$ a non negative, non
decreasing function such that $t^{-\alpha}\psi(t) \in
\mathcal{L}^\theta$. Then, for any $\delta
>0$ there exists a continuous and differentiable function $\varphi$ on
$\mathbb{R}_+$ such that:
\\i) $\psi(t) \le \varphi(t)$,
\\ii) $\varphi(t)t^{-\alpha -\delta}$ decreases and $\varphi(t)t^{-\alpha + \delta}$ increases,
\\iii) $\| t^{-\alpha}\varphi(t)\|_{\mathcal{L}^\theta} \le c \|t^{-\alpha} \psi(t) \|_{\mathcal{L}^\theta}$
where $c$ is a constant that only depends on $\delta$ and $\alpha$.
\end{lem}

\begin{proof} Follows the scheme of the one at \cite[Lemma
2.1]{KoPe}. We include it here only for completeness.

% We can supposse that
%$\delta < $.
 Set
\begin{equation*}
  \varphi_1(t) = (\alpha +\delta) t^{\alpha +\delta} \int_t^\infty u^{-\alpha -\delta} \psi(u)   \frac{du}{u}.
\end{equation*}
Then $\varphi_1(t)t^{-\alpha -\delta}$ decreases and
\begin{equation*}
  \varphi_1(t) \ge (\alpha + \delta) t^{\alpha +\delta} \psi(t) \int_t^\infty u^{-\alpha -\delta} \frac{du}{u} = \psi(t).
\end{equation*}
Furthermore, applying Hardy's inequality \cite[pg.124]{BeSh}, we
easily get that
\begin{equation}\label{l3}
  \|t^{-\alpha} \varphi_1(t)\|_{\mathcal{L}^\theta} \le c \| t^{-\alpha} \psi(t)   \|_{\mathcal{L}^\theta}.
\end{equation}
Set now
\begin{equation}\label{l4}
  \varphi(t) = 2\delta t^{\alpha -\delta} \int_0^t u^{-\alpha +\delta} \varphi_1(u) \frac{du}{u}
\end{equation}
Then $\varphi(t)t^{-\alpha + \delta}$ increases on $\mathbb{R}_+$ and
\begin{equation*}
  \varphi(t) \ge \varphi_1(t) \ge \psi(t) \quad t \in \mathbb{R}_+.
\end{equation*}
Furthermore, the change of variable $z= u^{2 \delta}$ in the right
hand side of (\ref{l4}) gives that
\begin{equation*}
  \varphi(t)t^{-\alpha -\delta} = t^{-2 \delta} \int_0^{t^{2\delta}} \mu(z^{1/(2 \delta)}) dz,
\end{equation*}
where  $\mu(u)=\varphi_1(u)u^{-\alpha - \delta}$ is a decreasing
function on $\mathbb{R}_+$. Thus, $\varphi(u)u^{-\alpha - \delta}$
decreases. Finally, using Hardy's inequality and (\ref{l3}) we get
(iii). The lemma is proved.
\end{proof}
Let $n \in \mathbb{N}$, $0 < r_j < \infty$, $1 \le p_j < \infty$ and $1 \le \theta_j \le \infty$ $\forall j \in \{ 1,\ldots,n\}$.
Denote
\begin{equation}\label{1}
r=n \left( \sum_{j=1}^n \frac{1}{r_j} \right)^{-1};\quad p=\frac{n}{r}\left( \sum_{j=1}^n \frac{1}{p_j r_j} \right)^{-1}; \quad \theta= \frac{n}{r} \left( \sum_{j=1}^n \frac{1}{\theta_j r_j} \right)^{-1}
\end{equation}
and
\begin{equation}\label{2}
  \beta_j= \frac{1}{r_j} \left( \frac{r}{n} + \frac{1}{p_j} - \frac{1}{p} \right).
\end{equation}
Then %$\beta_j < 1$ and
\begin{equation}\label{3}
  \sum_{j=1}^n \beta_j =1.
\end{equation}
%Indeed,
%\begin{equation*}
%  \left( \frac{r}{n} +\frac{1}{p_j} - \frac{1}{p} \right) \sum_{i=1}^n \frac{1}{r_i} = 1 + \sum_{i \neq j} \left( \frac{1}{p_j} - \frac{1}{p_i} \right) \frac{1}{r_i} < 1 + \sum_{i \neq j} \frac{1}{r_i} \le r_j \left( \sum_{i=1}^n \frac{1}{r_i} \right).
%\end{equation*}
%Thus, $\beta_j < 1$.
The equality (\ref{3}) follows immediately from
(\ref{1}).

In the following lemma we use the notations (\ref{1}) and (\ref{2}).
%%%%%%%%%%%%%%%%%%%%%%%%%%%%%%%%%%%%%%%%%%%%%%%%%
\begin{lem}\label{igualar}
%%%%%%%%%%%%%%%%%%%%%%%%%%%%%%%%%%%%%%%%%%%%%%%%%%
Let $n \in \mathbb{N}$, $0<r_j < \infty$, $1 \le p_j < \infty$ and
$1 \le \theta_j \le \infty$ for $j=1,\ldots,n$. Suppose that
$\beta_j>0$ for any $j$ and let
\begin{equation}\label{4}
  0 < \delta \le \frac{1}{2} \min_{1 \le j \le n} \{ \beta_j r_j \}.
\end{equation}
Let $\varphi_j$ positive, strictly increasing and continuously differentiable
 functions on $\mathbb{R}_+$, satisfying $\varphi_j(t)t^{-r_j} \in \mathcal{L}^{\theta_j}$. Besides $\varphi_j(t)t^{-r_j +\delta}$ increases and $\varphi_j(t)t^{-r_j-\delta}$ decreases.
Define
\begin{equation}\label{4bis}
  \sigma(t) = \inf \left\{ \max_{j=1,\ldots,n} \{t^{-1/p_j} \varphi_j (\delta_j) \} : \prod_{j=1}^n \delta_j =t,\,\, \delta_j >0 \right\}.
\end{equation}
Then:\\
i)There holds the inequality
\begin{equation}\label{5}
\left( \int_0^\infty t^{\theta(1/p -r/n)-1} \sigma(t)^\theta dt
\right)^{1/\theta} \le c \prod_{j=1}^n \left[ \|t^{-r_j}
\varphi_j(t) \|_{\mathcal{L}^{\theta_j}} \right]^{\frac{r}{n r_j}}.
\end{equation}
ii) There exist positive, continuously differentiable functions
$\delta_j(t)$ on $\mathbb{R}_+$ such that
\begin{equation}\label{6}
  \prod_{j=1}^n \delta_j(t) =t
\end{equation}
and
\begin{equation}\label{7}
  \sigma(t) = t^{-1/p_j} \varphi_j (\delta_j(t)) \quad (t\in \mathbb{R}_+,\, j=1,\ldots,n ).
\end{equation}
iii) for every $j=1,\ldots,n$
\begin{equation}\label{8}
 \sigma(t) t^{1/p -r/n +\delta} \uparrow \textnormal{ and } \sigma(t)t^{1/p -r/n -\delta} \downarrow ;
\end{equation}
\begin{equation}\label{9}
  \delta_j(t) t^{-\beta_j /3} \uparrow \textnormal{ and } \delta_j(t) t^{-3\beta_j} \downarrow ;
\end{equation}
iv) for every $j=1,\ldots,n$
\begin{equation}\label{10}
  \left( \int_0^\infty \left[ \frac{\varphi_j(\delta_j(t))}{\delta_j(t)^{r_j}}
  \right]^{\theta_j} \frac{dt}{t} \right)^{1/\theta_j} \le c \| t^{-r_j} \varphi_j(t)\|_{\mathcal{L}^{\theta_j}}.
\end{equation}
where c depends on $\delta$, $r_j$, $p_j$, $n$.
\end{lem}
%%%%%%%%%%%%%%%%%%%%%%%%%%%%%%%%%%%%%%%%%%%%%%%%%
% fin del lema de igualación
%%%%%%%%%%%%%%%%%%%%%%%%%%%%%%%%%%%%%%%%%%%%%%%%%

%%%%%%%%%%%%%%%%%%%%%%%%%%%%%%%%%%%%%%%%%%%%%%%%%
% demo del lema de igualaci\'on
%%%%%%%%%%%%%%%%%%%%%%%%%%%%%%%%%%%%%%%%%%%%%%%%%
\begin{proof} First note that
\begin{equation*}
\lim_{\delta_j \rightarrow 0} \varphi_j(\delta_j) =0 \quad \textnormal{and} \quad \lim_{\delta_j \rightarrow \infty} \varphi_j(\delta_j) = \infty.
\end{equation*}
Now we fix $t \in \mathbb{R}_+$. It's clear that there exists an unique point $\delta_j \equiv \delta_j(t) >0$ such that $\sigma(t) = t^{-1/p_j} \varphi_j(\delta_j(t))$.
Now let $1<\gamma < \infty$. Note that $\sigma(t) < t^{-1/p_j} \varphi_j (\delta_j(t) \gamma)$. Then
\begin{equation*}
  \sigma(t) < \min_{1 \le j \le n} \{ t^{-1/p_j} \varphi_j(\delta_j(t)\gamma) \}.
\end{equation*}
By the definition (\ref{4bis}) of $\sigma(t)$ we have that $\exists \delta_1^*,\ldots, \delta_n^*$ so that $\prod_{j=1}^n \delta_j^* =t$ and
\begin{equation*}
  \sigma(t) \le \max_{j=1,\ldots,n} \{ t^{-1/p_j} \varphi_j(\delta_j^*) \} < \min_{1 \le j \le n} \{ t^{-1/p_j} \varphi_j(\delta_j(t)\gamma) \}.
\end{equation*}
Therefore
\begin{equation*}
\delta_j^* < \delta_j(t) \gamma \quad \textnormal{for all } \,
j\quad \Longrightarrow t = \prod_{j=1}^n \delta_j^* < \gamma^n
\prod_{j=1}^n \delta_j(t).
\end{equation*}
Taking limits when $\gamma$ tends to $1$ we have $t \le \prod_{j=1}^n \delta_j(t)$.
Now if $\prod_{j=1}^n \delta_j(t) >t$ we choose $0< \delta'_j < \delta_j(t)$ so that $\prod_{j=1}^n \delta'_j =t$. We have $\sigma(t) = t^{-1/p_j} \varphi_j(\delta_j(t)) > t^{-1/p_j} \varphi_j(\delta'_j)$. Then
\begin{equation*}
  \sigma(t) > \max_{j=1,\ldots,n} \{t^{-1/p_j} \varphi_j(\delta'_j) \} \textnormal{ and } \prod_{j=1}^n \delta'_j =t
\end{equation*}
which contradicts (\ref{4bis}) definition of $\sigma(t)$.
So functions $\delta_j(t)$ satisfy (\ref{6}) and (\ref{7}).

Besides, for any $j=1,\ldots,n$, by (\ref{7})
\begin{equation}\label{11}
  \delta_j(t) = \varphi_j^{-1}( t^{1/p_j - 1/p_n} \varphi_n(\delta_n(t))).
\end{equation}
Then by (\ref{6})
\begin{equation*}
t = \Phi(\delta_n(t),t)
\end{equation*}
where
\begin{equation*}
  \Phi(s,t)= s \prod_{j=1}^{n-1} \varphi_j^{-1}( t^{1/p_j - 1/p_n} \varphi_n(s))
\end{equation*}
which is a function of $C^1(\mathbb{R}_+^2)$ strictly increasing
respect to $s$. In consequence, by the implicit function theorem we
have that
 $\delta_n \in C^1(\mathbb{R}_+)$ and so, by (\ref{11}) $\delta_j \in C^1(\mathbb{R}_+)$ for all $j =1,\ldots,n$.
We have just proved (ii).

Our conditions on $\varphi_j$ implies that for every $j=1,\ldots,n$
\begin{equation}\label{12}
 \frac{-r_j -\delta}{t} \le - \frac{\varphi_j'(t)}{\varphi_j(t)} \le \frac{-r_j +\delta}{t}.
\end{equation}
Besides
\begin{equation}\label{13}
 \frac{\sigma'(t)}{\sigma(t)} = \frac{-1/p_j}{t} + \frac{\varphi_j'(\delta_j(t))}{\varphi_j(\delta_j(t))} \delta_j'(t).
\end{equation}
Now we derive (\ref{6}) and taking into account (\ref{3}) we have that for every $t>0$ $\exists$ $m \equiv m(t)$ and $l \equiv l(t)$ such that
\begin{equation*}
 \frac{\delta_m'(t)}{\delta_m(t)} \le \frac{\beta_m}{t} \textnormal{ and } \frac{\delta_l'(t)}{\delta_l(t)} \ge \frac{\beta_l}{t}.
\end{equation*}
Then
\begin{equation*}
  \frac{-1/p_l}{t} + \frac{\varphi_l'(\delta_l(t))}{\varphi_l(\delta_l(t))} \delta_l(t) \frac{\beta_l}{t} \le \frac{\sigma'(t)}{\sigma(t)} \le \frac{-1/p_m}{t} + \frac{\varphi_m'(\delta_m(t))}{\varphi_m(\delta_m(t))} \delta_m(t) \frac{\beta_m}{t}.
\end{equation*}
Now, using (\ref{12}) we obtain
\begin{equation*}
\frac{-1/p_l + \beta_l(r_l - \delta)}{t} \le \frac{\sigma'(t)}{\sigma(t)} \le \frac{-1/p_m +\beta_m(r_m+ \delta)}{t}.
\end{equation*}
And due to $0< \beta_l,\beta_m < 1$ and (\ref{2}) we obtain:
\begin{equation}\label{14}
  \frac{r/n -1/p -\delta}{t} \le \frac{\sigma'(t)}{\sigma(t)} \le \frac{r/n -1/p
  +\delta}{t},
\end{equation}
which implies (\ref{8}). Besides, for (\ref{14}) and (\ref{13})
\begin{equation*}
  \frac{\beta_j r_j - \delta}{t} \le \frac{\varphi_j'(\delta_j(t))}{\varphi_j(\delta_j(t))}
  \delta_j'(t) \le \frac{\beta_j r_j +\delta}{t}.
\end{equation*}
Then $\delta_j$ increases and using this last inequality and
(\ref{12})
\begin{equation}\label{15}
   \frac{\beta_j r_j - \delta}{(r_j +\delta)t}
    \le \frac{\delta_j'(t)}{\delta_j(t)} \le \frac{\beta_j r_j + \delta}{(r_j - \delta)t}.
\end{equation}
From here and (\ref{4})
\begin{equation*}
  \frac{\beta_j}{3t} \le \frac{\delta_j'(t)}{\delta_j(t)} \le \frac{3
  \beta_j}{t};
\end{equation*}
and we have proved (\ref{9}). The statement (iv) is the immediate
consequence of applying the left inequality of (\ref{15}) and the
change of variable $u=\delta_j(t)$.

Finally we multiply (\ref{7}) elevated to $1/r_j$ and use (\ref{6})
\begin{equation*}
  \sigma(t)^{n/r} = t^{-n/(rp)} \prod_{j=1}^n \varphi_j(\delta_j(t))^{1/r_j} = t^{-n/(rp) +1} \prod_{j=1}^n \left[ \frac{\varphi_j(\delta_j(t))}{\delta_j(t)^{r_j}} \right]^{1/r_j}.
\end{equation*}
So
\begin{equation*}
 \left( \int_0^\infty t^{\theta (1/p -r/n) -1} \sigma(t)^\theta dt \right)^{1/\theta} =
  \left( \int_0^\infty \prod_{j=1}^n \left[ \frac{\varphi_j(\delta_j(t))}{\delta_j(t)^{r_j}} \right]^\frac{\theta r}{n r_j} \frac{dt}{t} \right)^{1/\theta} \le
\end{equation*}
\begin{equation*}
\le  \prod_{j=1}^n \left( \int_0^\infty \left[
\frac{\varphi_j(\delta_j(t))}{\delta_j(t)^{r_j}} \right]^{\theta_j}
\frac{dt}{t} \right)^{\frac{1}{\theta_j} \frac{r}{nr_j}}
\end{equation*}
and applying (\ref{10}) we obtain (\ref{5}).
The lemma is proved.
\end{proof}
%%%%%%%%%%%%%%%%%%%%%%%%%%%%%%%%%%%%%%%%%%%
% fin de la demo del lema de igualaci\'on
%%%%%%%%%%%%%%%%%%%%%%%%%%%%%%%%%%%%%%%%%%%

In the following lemma we use the notations (\ref{1}) and (\ref{2}) too.
%%%%%%%%%%%%%%%%%%%%%%%%%%%%%%%%%%%%%%%%%%%
% lema tipo de inmersión con exponente límite
%%%%%%%%%%%%%%%%%%%%%%%%%%%%%%%%%%%%%%%%%%%%
\begin{lem}\label{explim}
Let $n \in \mathbb{N}$, $0< r_j < \infty$, $1 \le p_j < \infty$ and
$1 \le \theta_j \le \infty$ $(j=1,\ldots,n) $ such that $\beta_j
>0$ for all $j$. Then for every function $f \in
S_0(\mathbb{R}^n)\cap L_{loc}(\mathbb{R}^n)$ such that $f \in
b_{p_1,\ldots,p_n; \theta_1, \ldots,
\theta_n}^{r_1,\ldots,r_n}(\mathbb{R}^n)$
\begin{equation}\label{16}
  f^*(t) \le (2^{k'}-1)f^*(\xi t) + c(\xi) \sigma(t), \quad \textnormal{for all } \xi >1
\end{equation}
and
\begin{equation}\label{17}
  \left( \int_0^\infty t^{\theta (1/p -r/n) -1} \sigma(t)^\theta dt \right)^{1/\theta} \le \prod_{j=1}^n \|h^{-r_j} w_j^{k_j}
  (f;h)_{p_j}\|_{\mathcal{L}^{\theta_j}}^
  \frac{r}{nr_j},
\end{equation}
where $r_j < k_j \in \mathbb{N}$, $k'=\max k_j $.
\end{lem}
%%%%%%%%%%%%%%%%%%%%%%%%%%%%%%%%%%%%%%%%%%%%
% demo del lema tipo exponente limite
%%%%%%%%%%%%%%%%%%%%%%%%%%%%%%%%%%%%%%%%%%%%
\begin{proof} First we apply Lemma \ref{eastj} to $f$ getting
its estimation (\ref{l2}). Now we define $0< \delta = \frac{1}{2}
\min_j \{ \beta_j r_j \}$ and apply Lemma \ref{crec} to the modulus
of continuity getting
\begin{equation*}
  w_j^{k_j}(f;u)_{p_j} \le \varphi_j(u), \quad \varphi_j(u)u^{-r_j + \delta} \uparrow , \quad \varphi_j(u)u^{-r_j -\delta} \downarrow
\end{equation*}
and
\begin{equation}\label{18}
 \| \varphi_j(u)u^{-r_j}\|_{\mathcal{L}^{\theta_j}} \le c \|h^{-r_j} w_j^{k_j}(f;h)_{p_j}\|_{\mathcal{L}^{\theta_j}}.
\end{equation}
Then we obtain (\ref{16}) with $\sigma(t)$ of (\ref{4bis}) kind.
Only rests to apply Lemma \ref{igualar} and use (\ref{18}).
\end{proof}

\section{Results}\label{secteo}

%%%%%%%%%%%%%%%%%%%%%%%%%%%%%%%%%%%%%%%%%
% inmersión sin exponente límite
%%%%%%%%%%%%%%%%%%%%%%%%%%%%%%%%%%%%%%%%%
\begin{teo}[embedding with no limit exponent]\label{nolimit}
Assume that a function $f$ satisfies the conditions of Lemma \ref{explim} and $f \in L^1(\mathbb{R}^n) + L^{p_0}(\mathbb{R}^n)$. For some $p_0>0$ such that
\begin{equation*}
 \frac{1}{p_0} > \frac{1}{p} - \frac{r}{n}.
\end{equation*}
Let $\max(1,p_0)< q < \infty$ and
\begin{equation}\label{234}
 \frac{1}{q} > \frac{1}{p} - \frac{r}{n}.
\end{equation}
Then for any $s >0$ $f \in L^{q,s}(\mathbb{R}^n)$ and
\begin{equation}\label{235}
 \|f \|_{q,s} \le c \left[ \|f\|_{L^1 +L^{p_0}} + \prod_{j=1}^n \|f\|_{b_{p_j,\theta_j;j}^{r_j}} \right].
\end{equation}
\end{teo}
%%%%%%%%%%%%%%%%%%%%%%%%%%%%%%%%%%%%%%%
% fin del corolario
%%%%%%%%%%%%%%%%%%%%%%%%%%%%%%%%%%%%%%%
\begin{proof} The proof follows the scheme of the one of \cite[Corollary 2.5]{KoPe}.
 We include it here for the reader's convenience. We can assume that $s< \min(1,p_0,\theta)$. Let $f=g+h$, with
$g \in L^{1}(\mathbb{R}^n)$ and $h \in L^{p_0}(\mathbb{R}^n)$.
Applying the H\"older inequality, we obtain
\begin{equation*} J_1 \equiv \int_1^\infty \left[ t^{1/q} f^*(t)
\right]^s
 \frac{dt}{t} \le c \, \left[ \left( \int_0^\infty g^*(t) dt \right)^{s}
+ \left(\int_0^\infty h^*(t)^{p_0} dt \right)^{s /p_0}
\right].\end{equation*} It follows that

\begin{equation}\label{236} J_1 \le c^\prime \, \|f\|_{L^1 + L^{p_0}}^s .
\end{equation} Let $0<\delta<1$. Using (\ref{16}) with
$\xi=(2^{1/s}K)^q$, we get by H\"older inequality and (\ref{234}):
\begin{equation*} J_\delta \equiv \int_\delta^\infty \left[ t^{1/q} f^*(t) \right]^s
 \frac{dt}{t} \le J_1 + K^s \int_\delta^1 \left[ t^{1/q} f^*(\xi t) \right]^s \frac{dt}{t} + \end{equation*}
\begin{equation*} + c\int_0^1 t^{s/q \, -1} \sigma(t)^s dt \le J_1+ \frac{1}{2}J_\delta + \end{equation*}
\begin{equation*}+ c^\prime \left(\int_0^1 t^{\theta(1/p \, -r/n)-1} \sigma(t)^\theta dt \right)^{s/\theta}.\end{equation*}
The inequality (\ref{235}) follows now from (\ref{17}) and (\ref{236}).
\end{proof}
%%%%%%%%%%%%%%%%%%%%%%%%%%%%%%%%%%%%%%%%
% fin de la demo del corolario
%%%%%%%%%%%%%%%%%%%%%%%%%%%%%%%%%%%%%%%%

\begin{teo}[embedding with limit exponent]\label{limit}
Let $n \in \mathbb{N}$, $0< r_j < \infty$, $1 \le p_j < \infty$, $1 \le \theta_j \le \infty$.
 Define $r$, $p$ and $\theta$ as in (\ref{1}) and supposse that (see (\ref{2})) $\beta_j > 0$
  $(1\le j\le n)$ and $p< \frac{n}{r}$. We define $q^*=np/(n-rp)$. Then, for every function
   $f \in L_{loc}(\mathbb{R}^n)$ there holds
\begin{equation*}
 \| f \|_{q^*,\theta} \le c \prod_{j=1}^n
 \left[\|f\|_{b_{p_j,\theta_j;j}^{r_j}}\right]^\frac{r}{nr_j}
 %\prod_{j=1}^n \left( \int_0^\infty [h^{-r_j} w_j^{k_j}(f;h)_{p_j}]^{\theta_j} \frac{dh}{h} \right)^{\frac{1}{\theta_j} \frac{r}{nr_j}}
\end{equation*}
\end{teo}
This statement follows immediately from Lemma \ref{explim}. Let's
note that the condition of $f\in L_{loc}(\mathbb{R}^n)$ and being of
compact support can be substituted by $f\in L^{p_0}(\mathbb{R}^n)$
for some $0<p_0<q^*.$ In this case we have the embedding
\begin{equation*}
  L^{p_0}\cap b_{p_1,\ldots,p_n;\theta_1,\ldots,\theta_n}^{r_1,\ldots,r_n}(\mathbb{R}^n)\hookrightarrow
  L^{q^*,\theta}(\mathbb{R}^n).
\end{equation*}

The following theorems is our main result. It expresses an embedding
of different metrics.
\begin{teo}\label{maint}
Let $n \in \mathbb{N}$, $0< r_j < \infty$, $1 \le p_j < \infty$, $1
\le \theta_j \le \infty$ $j \in \{1,\ldots,n\}$. Let $r$, $p$ and
$\theta$ be the numbers defined in (\ref{1}) and supposse that for
all $1\le j\le n$
\begin{equation*}
  \beta_j = \frac{1}{r_j} \left( \frac{r}{n} +\frac{1}{p_j} -\frac{1}{p} \right) >0.
\end{equation*}
We choose arbitrary $p_j < q_j < \infty$ such that
\begin{equation*}
  \frac{1}{q_j} > \frac{1}{p} - \frac{r}{n}
\end{equation*}
and denote
\begin{equation*}
 \varkappa_j = 1- \frac{1}{\beta_j r_j} \left( \frac{1}{p_j} - \frac{1}{q_j} \right),
\end{equation*}
\begin{equation*}
  \alpha_j = \varkappa_j r_j \quad \textnormal{and} \quad \frac{1}{\theta'_j} = \frac{1- \varkappa_j}{\theta} + \frac{\varkappa_j}{\theta_j}.
\end{equation*}
Then for any function $f \in S_0(\mathbb{R}^n)$ such that $f \in
b_{p_1,\ldots,p_n;\theta_1,\ldots,\theta_n}^{r_1,\ldots,r_n}
(\mathbb{R}^n)$ %$f \in B_{p_i,\theta_i;i}^{r_i}(\mathbb{R}^n)\,\,
%\forall i=1,\ldots,n.$
 there holds the inequality
\begin{equation}\label{19pre}
 \left( \int_0^\infty [ h^{-\alpha_j} \| \Delta_j^{r_j} (h)f \|_{q_j,1}]^{\theta'_j} \frac{dh}{h} \right)^{1/\theta'_j}
 \le c \sum_{i=1}^n %\left( \int_0^\infty [h^{-r_i} w_i^{k_i}(f;h)_{p_i}]^{\theta_i} \frac{dh}{h} \right)^{/\theta_i}
\|f\|_{b_{p_i,\theta_i;i}^{r_i}}
\end{equation}
(where $r_j<k_j\in\mathbb{N}$); which implies the embedding
\begin{equation}\label{desmetricas}
b_{p_1,\ldots,p_n;\theta_1,\ldots,\theta_n}^{r_1,\ldots,r_n}(\mathbb{R}^n)
\hookrightarrow b_{q_1,\ldots,q_n;\theta'_1, \ldots,
\theta'_n}^{\alpha_1,\ldots,\alpha_n}(\mathbb{R}^n).
\end{equation}
\end{teo}
%%%%%%%%%%%%%%%%%%%%%%%%%%%%%%%%%%%%%%%%%%%%
% fin del enunciado del teorema
%%%%%%%%%%%%%%%%%%%%%%%%%%%%%%%%%%%%%%%%%%%%
\begin{proof}
Supposse that $j=1$. Note that $0< \varkappa_1 < 1$.
Set now $r_i < k_i \in \mathbb{N}$ for every $i=1,\ldots,n$.
Denote for $h>0$ and $x\in \mathbb{R}^n$
\begin{equation*}
 f_h(x) = | \Delta_1^{k_1}(h)f(x)|.
\end{equation*}
As $\|f_h\|_{p_1} \le w_1^{k_1}(f;h)_{p_1} < + \infty$, $f_h \in
L^{p_1}(\mathbb{R}^n)$ and applying Theorem \ref{nolimit} we have
that $f_h \in L^{q_1,1}(\mathbb{R}^n)$. Denote for $h>0$
\begin{equation*}
  J(h) \equiv \|f_h \|_{q_1,1} = \int_0^\infty t^{1/q_1 -1} f_h^*(t) dt<\infty.
\end{equation*}
Set $\xi_0 = (2^{k+2})^{q_1}$ and
\begin{equation}\label{19}
  Q(h) = \{ t>0: f_h^*(t) \ge 2^{k+1}f_h^*(\xi_0 t) \},
\end{equation}
where $k=\max k_i$ as in Lemma \ref{eastj}. Then
\begin{equation*}
  \int_{\mathbb{R}_+ - Q(h)} t^{1/q_1 -1} f_h^*(t) dt \le 2^{k+1} \int_0^\infty t^{1/q_1 -1} f_h^*(\xi_0 t)dt =
\end{equation*}
\begin{equation*}
 = 2^{k+1} \xi_0^{-1/q_1} \int_0^\infty t^{1/q_1 -1} f_h^*(t)dt = \frac{1}{2} J(h).
\end{equation*}
Therefore
\begin{equation*}
  J(h) \le 2 \int_{Q(h)} t^{1/q_1 -1} f_h^*(t) dt \equiv 2 J'(h).
\end{equation*}
So, it is enough to estimate $f_h$ in $Q(h)$. Now we choose
\begin{equation}\label{20}
  0< \delta = \frac{1}{2} \min \{ \min_i \{ \beta_i r_i \}, \frac{1}{q_1} - \frac{1}{p} +\frac{r}{n} \}
\end{equation}
By virtue of Lemma \ref{crec} there exists functions $\varphi_i(t)$ on $\mathbb{R}_+$ $(i=1,\ldots,n)$
 continuously differentiable such that
\begin{equation}\label{21}
 \varphi_i(t)t^{-r_i - \delta} \downarrow \textnormal{ and } \varphi_i(t)t^{-r_i +\delta}
 \uparrow;
\end{equation}
\begin{equation}\label{22}
  w_i^{k_i}(f;t)_{p_i} \le \varphi_i(t);
\end{equation}
\begin{equation}\label{22bis}
  \| t^{-r_i} \varphi_i(t) \|_{\mathcal{L}^{\theta_i}} \le c \|t^{-r_i} w_i^{k_i}(f;t)_{p_i}\|_{\mathcal{L}^\theta_i}=c\|f\|_{b^{r_i}_{p_i,\theta_i;i}}.
\end{equation}
Now, applying Lemma \ref{eastj}, (\ref{22}) and (\ref{19}) we have
for all $t \in Q(h)$
\begin{equation*}
  f_h^*(t) \le c \max_{1 \le i \le n} t^{-1/p_i} \varphi_i(\delta_i(t)),
\end{equation*}
with $\delta_i$ any functions on $\mathbb{R}_+$ such that $
\prod_{i=1}^n \delta_i(t) = t$.

 Due to (\ref{20}) and (\ref{21}) we apply Lemma \ref{igualar}
and we have that there exists a non-negative function $\sigma(t)$
such that
\begin{equation}\label{23}
  f_h^*(t) \le c \sigma(t)
\end{equation}
\begin{equation}\label{24}
  \left( \int_0^\infty t^{\theta (1/p -r/n)-1} \sigma(t)^\theta dt \right)^{1/\theta} \le
   \prod_{i=1}^n \left[ \| \varphi_i (t) t^{-r_i} \|_{\mathcal{L}^{\theta_i}}
   \right]^\frac{r}{nr_i}.
\end{equation}
There exist positive, continuously differentiable functions $u_i(t)$
on $\mathbb{R}_+$ such that $\prod_{i=1}^n u_i(t) =t$ and
\begin{equation}\label{25}
  \sigma(t) = t^{-1/p_i} \varphi_i(u_i(t)) \quad \forall i=1,\ldots,n.
\end{equation}
\begin{equation}\label{25bis}
  \sigma(t) t^{1/p -r/n +\delta} \uparrow,
\end{equation}
\begin{equation}\label{26}
 u_1(t)t^{-\frac{\beta_1}{3}} \uparrow \quad \textnormal{and} \quad u_1(t)t^{-3\beta_1}
 \downarrow,
\end{equation}
\begin{equation}\label{27}
 \left( \int_0^\infty [u_1(t)^{-r_1} \varphi_1(u_1(t)) ]^{\theta_1} \frac{dt}{t} \right)^{1/\theta_1} \le
  c \| \varphi_1(t)t^{-r_1}\|_{\mathcal{L}^{\theta_1}}.
\end{equation}
The estimation (\ref{23}) can be used for``little" $t$. For ``big"
$t$ we will use the following estimate, which is consequence of a
weak type inequality and (\ref{22}).
\begin{equation}\label{28}
 f_h^*(t) \le t^{-1/p_1} \| f_h \|_{p_1} \le t^{-1/p_1} w_1^{k_1}(f;h)_{p_1} \le t^{-1/p_1} \varphi_1(h).
\end{equation}
Then
\begin{equation}\label{29}
  J'(h) \le \int_0^\infty t^{1/q_1 -1} \Phi(t,h) dt
\end{equation}
where
\begin{equation}\label{30}
  \Phi(t,h) = \min \{ \sigma(t), t^{-1/p_1} \varphi_1(h) \}.
\end{equation}
Then
\begin{multline*}
  J \equiv \|h^{-\alpha_1} J(h)\|_{\mathcal{L}^{\theta'_1}} \le\\c\left(\| h^{-\alpha_1 } \int_{\{h \le u_1(t)\}  } t^{1/q_1 -1/p_1 -1}
   \varphi_1(h) dt \|_{\mathcal{L}^{\theta'_1}}\right.   +\\+\left.
 \| h^{-\alpha_1 }  \int_{\{ h \ge u_1(t)\} } t^{1/q_1
-1} \sigma(t) dt \|_{\mathcal{L}^{\theta'_1}}\right)  \equiv c[ J_1
+J_2].
\end{multline*}

 Due to (\ref{26}) $u_1$ possess positive inverse $u_1^{-1}
\equiv \beta$ on $\mathbb{R}_+$ and
\begin{equation}\label{31}
  \frac{u_1'(t)}{u_1(t)} \le \frac{c}{t}.
\end{equation}
Now
\begin{equation}\label{32}
  J_1 = \| h^{-\alpha_1} \varphi_1(h)
   \int_{\beta(h)}^\infty t^{1/q_1 -1/p_1 -1} dt \|_{\mathcal{L}^{\theta'_1}}
= c \| h^{-\alpha_1} \varphi_1(h) \beta(h)^{ (1/q_1
-1/p_1)}\|_{\mathcal{L}^{\theta'_1}}.
\end{equation}

 Now we proceed with $J_2$. Due to (\ref{25bis}) and
(\ref{20})
\begin{multline*}
 J_2 = \| h^{-\alpha_1} \int_0^{\beta(h)} t^{1/q_1-1} \sigma(t) dt \|_{\mathcal{L}^{\theta'_1}}  \le
\\
  \le \| h^{-\alpha_1} \beta(h)^{1/p-r/n +\delta} \sigma(\beta(h)) \int_0^{\beta(h)} t^{1/q_1 -1/p +r/n -\delta-1} dt \|_{\mathcal{L}^{\theta'_1}} =
\\
  = c \| h^{-\alpha_1} \beta(h)^{1/q_1} \sigma(\beta(h))\|_{\mathcal{L}^{\theta'_1}}
\end{multline*}

By (\ref{25}) the last integral is the same as one in the right hand of (\ref{32}).
Therefore we have that
\begin{equation*}
  J \le c\| h^{-\alpha_1} \beta(h)^{ \frac{1}{q_1} - \frac{\varkappa_1}{p_1}  } \varphi_1(h)^{\varkappa_1} \sigma( \beta(h))^{(1-\varkappa_1)}\|_{\mathcal{L}^{\theta'_1}}.
\end{equation*}
The change of variable $u_1(z) =h$ and (\ref{31}) arrive us to
\begin{equation*}
 J \le c \|u_1(z)^{-\alpha_1} \varphi_1( u_1(z))^{\varkappa_1 } z^{  \frac{1}{q_1} - \frac{\varkappa_1 }{p_1} }
 \sigma(z)^{(1-\varkappa_1)}\|_{\mathcal{L}^{\theta'_1}}.
\end{equation*}
And using H\"older's inequality with exponents $u= \theta_1/(\varkappa_1 \theta'_1)$ and $u'= \theta_1 /(\theta_1 - \varkappa_1 \theta'_1)$ (observe that $(1-\varkappa_1)\theta'_1 u' = \theta$, $(\frac{\theta'_1}{q_1} - \frac{\varkappa_1 \theta'_1}{p_1})u' = \theta ( \frac{1}{p} - \frac{r}{n})$).
\begin{equation*}
  J \le \left( \int_0^\infty \left[ \frac{\varphi_1(u_1(z))}{u_1(z)^{r_1}} \right]^{\theta_1} \frac{dz}{z} \right)^{\varkappa_1 / \theta_1} \left( \int_0^\infty z^{\theta (1/p-r/n)} \sigma(z)^\theta \frac{dz}{z}
  \right)^{(1-\varkappa_1)/\theta}.
\end{equation*}
From here and (\ref{27}), (\ref{24}) and (\ref{22bis}) follows
\begin{equation*}
  J \le \left( \|t^{-r_1} w_1^{k_1}(f;t)_{p_1} \|_{\mathcal{L}^{\theta_1}}\right)^{\varkappa_1}
   \prod_{i=1}^n \left( \|t^{-r_i} w_i^{k_i}(f;t)_{p_i}\|_{\mathcal{L}^{\theta_i}} \right)^{\frac{r(1-\varkappa_1)}{n
   r_i}}.
\end{equation*}
which, using the inequality between arithmetic and geometric means,
implies (\ref{19pre}). The theorem is proved.
\end{proof}

\begin{rem}\label{mejorlorentz}
  Let us note that in (\ref{19pre}) it appears the stronger Lorentz
  norm $L^{q_j,1}$ instead of the norm $L^{q_j}$. Then
  (\ref{19pre}) is stronger than (\ref{desmetricas}). A detailed
  reading of the proof shows that in fact, in the right part of
  (\ref{19pre}) can appear the $L^{q_j,\xi}$ norm, for every
  $\xi>0$. Note that in this case ``$c$" in the right part of
  (\ref{19pre}) explodes when $\xi$ goes to $0$.
\end{rem}

\begin{rem}
  The values of the parameters $\theta'_j$ found in Theorem \ref{maint} are sharp.
   In order to see that the values of $\theta'_j$ can not be smaller
   one can consider the close embedding relation between Sobolev and
   Besov spaces and the fact that the parameters found in
   \cite{KoPe} are sharp (\cite[Remark 3.3]{KoPe}).
\end{rem}

%%%%%%%%%%%%%%%%%%%%%%%%%%%%%%%%%%%%%%%%%%%%%%%%%%

\begin{thebibliography}{99}

\bibitem {BeSh} {\sc Bennett, C.}, and {\sc Sharpley, R.}:
 Interpolation of Operators, Academic Press, 1988.

\bibitem {BeIlNi} {\sc Besov, O.V., Il'in, V.P.}, and {\sc Nikol'ski\u\i, S.M.}:
Integral Representation of Functions and Imbedding Theorems, vol. 1
-- 2,
 Winston, Washington D.C., Halsted, New York--Toronto--London, 1978

\bibitem{Go}{\sc Golovkin, K.K.}: A generalization of
Marcinkiewicz's interpolation theorem, Trudy Mat. Inst. Steklov {\bf
102} (1967), 5 -- 28.

\bibitem    {He} {\sc Herz, C.}: Lipschitz spaces and Bernstein's
theorem of absolutely convergent Fourier transform, J. Math. Mech.
{\bf 18} No. 18 (1968), 283 -- 323.

\bibitem {Ko93} {\sc Kolyada, V.I.}: On embedding of Sobolev spaces, Mat. Zametki {\bf 54} No. 3 (1993), 48 -- 71; English transl. in Math. Notes {\bf 54} No. 3 (1993), 908 -- 922


\bibitem {Ko98} {\sc Kolyada, V.I.}: Rearrangements of functions
 and embedding of anisotropic spaces of the Sobolev type, East J. on
  Approximations {\bf 4} No. 2 (1998), 111 -- 199.

\bibitem {Ko01} {\sc Kolyada, V.I.}: Embeddings of fractional Sobolev spaces and estimates of Fourier transforms, Mat. Sb. {\bf 192} No. 7 (2001), 51 -- 72; English transl. in Sbornik: Mathematics {\bf 192} No. 7 (2001), 979 -- 1000


\bibitem {KoPe} {\sc Kolyada, V.I., and Pérez, F.J.}: Estimates of
difference norms for functions in anisotropic Sobolev spaces, Mat.
Nachr. {\bf 267}, 46--64 (2004).

\bibitem {Nik51} {\sc Nikol'ski\u\i, S.M.}: Inequalities for
entire functions of finite degree and their application in the
theory of differentiable functions of several variables, Trudy Mat.
Inst. Steklov {\bf 38} (1951), 244 -- 278 (in Russian).

\bibitem {Ni} {\sc Nikol'ski\u\i, S.M.}: Approximation of Functions
 of Several Variables and Imbedding Theorems, Springer -- Verlag, Berlin -- Heidelberg -- New York, 1975

\bibitem {Pee} {\sc Peetre, J.}: Espaces d'interpolation et espaces de Soboleff, Ann. Inst. Fourier (Grenoble) {\bf 16} (1966), 279 -- 317

\bibitem{Pe} {\sc P\'erez L\'azaro, F.J.}: A note on extreme cases
of Sobolev embeddings, J. Math. Anal. Appl. 320 (2006) 973 -- 982.

\bibitem{Ul} {\sc Ul'yanov, P.L.}: On the embedding of certain classes of functions, Mat. Zametki {\bf 1} 4(1967), 405--414.


\end{thebibliography}
\end{document}